 \documentclass[12pt]{article}
 \usepackage{amsmath,amsfonts,theorem}
 \numberwithin{equation}{section}

 \newtheorem{prop}{Proposition}[section]

 \newtheorem{cor}{Corollary}[section]

 \newtheorem{thm}{Theorem}[section]

 \theorembodyfont{\rmfamily}
 \newtheorem{dfn}{Definition}[section]

 \newtheorem{rmk}{Remark}

 \newcommand{\qed}{\ifhmode\unskip\nobreak\fi\quad\ensuremath\square}

 \newcommand{\Span}[1]{\left< #1 \right>}
 
 \newcommand{\dd}{\mathrm{d}}
 \newcommand{\id}{\operatorname{id}}
 \newcommand{\pr}{\operatorname{pr}}

 \newcommand{\CS}{\operatorname{CS}}
 \newcommand{\BS}{\mathrm{BS}}
 \newcommand{\kBS}{k\mathrm{\text{-}BS}}
 \newcommand{\HP}{\mathrm{HP}}
 
 \newcommand{\LC}{\mathrm{LC}}
 \newcommand{\US}{\mathrm{US}}

 \newcommand{\semis}{^{\mathrm{ss}}}

 \newcommand{\irr}{^{\mathrm{irr}}}
 \newcommand{\red}{^{\mathrm{red}}}
 \newcommand{\rest}[1]{{}_{{\textstyle{|}}#1}} 

 \newcommand{\albar}{\overline{\al}}
 \newcommand{\dbar}{\overline{\partial}}
 \newcommand{\Dbar}{\overline{D}}
 \newcommand{\gbar}{\overline{g}}

 \newcommand{\p}{\partial}
 
 \newcommand{\tensor}{\otimes}

 \newcommand{\sA}{\mathcal A} 
 \newcommand{\sE}{\mathcal E}
 \newcommand{\sG}{\mathcal G} 
 \newcommand{\sF}{\mathcal F}
 \newcommand{\sH}{\mathcal H}
 \newcommand{\sL}{\mathcal L}
 \newcommand{\sM}{\mathcal M}
 \newcommand{\Oh}{\mathcal O}
 \newcommand{\sS}{\mathcal S}
 \newcommand{\sW}{\mathcal W}

 \newcommand{\al}{\alpha}
 
 \newcommand{\de}{\delta}
 \newcommand{\ep}{\varepsilon}
 \newcommand{\fie}{\varphi}

 \newcommand{\om}{\omega}

 \newcommand{\Ga}{\Gamma}
 \newcommand{\La}{\Lambda}
 \newcommand{\Om}{\Omega}
 
 \newcommand{\Si}{\Sigma}
 \newcommand{\na}{\nabla}

 \newcommand{\PP}{\mathbb P}
 \newcommand{\C}{\mathbb C}

 \newcommand{\R}{\mathbb R}
 \newcommand{\Z}{\mathbb Z}

 \newcommand{\diag}{\operatorname{diag}}
 \newcommand{\rk}{\operatorname{rank}}

 \newcommand{\Diff}{\operatorname{Diff}}
 \newcommand{\tr}{\operatorname{tr}}
 \newcommand{\End}{\operatorname{End}}
 \newcommand{\GFT}{\operatorname{GFT}}
 \newcommand{\sGFT}{\operatorname{sGFT}}

 \newcommand{\Hol}{\operatorname{Hol}}
 \newcommand{\Hom}{\operatorname{Hom}}
 \newcommand{\Lie}{\operatorname{Lie}}
 \newcommand{\Pic}{\operatorname{Pic}} 

 \newcommand{\Sing}{\operatorname{Sing}}
 \newcommand{\Vol}{\operatorname{Vol}}

 \newcommand{\U}{\operatorname U} 
 \newcommand{\PU}{\operatorname{PU}}
 \newcommand{\SU}{\operatorname{SU}}
 \newcommand{\fsl}{\operatorname{\mathfrak{su}}}
 \newcommand{\su}{\operatorname{\mathfrak{su}}}

\begin{document}

 \title{Quantization and ``theta functions''}
 \markright{\hfill Quantization and ``theta functions'' \quad}

 \author{Andrei Tyurin}
 \date{Apr 1999}
 \maketitle

 \begin{abstract}
 Geometric Quantization links holomorphic geometry with real geo\-metry, a
relation that is a prototype for the modern development of mirror
symmetry. We show how this treatment can be used to construct a special
basis in every space of conformal blocks. This is a direct generalization
of the basis of theta functions with characteristics in every complete
linear system on an Abelian variety (see \cite{Mum}). The same construction
generalizes the classical theory of theta functions to vector bundles of
higher rank on Abelian varieties and K3 surfaces. We also discuss the
geometry behind these constructions.
 \end{abstract}

 \section{Introduction}
 \markright{\hfill Quantization and ``theta functions'' \quad}

It is a fruitful question to ask for some special basis of the complete
linear systems $\PP H^0(X, L^k)$, where $X$ is a smooth complete algebraic
variety and $L$ a polarization. After this, following Mumford, we can ask
for special equations defining $X$ under its embedding in $\PP
H^0(X,L^k)^*$. Of course, this is a priori impossible (for example, for
$\PP H^0(\PP^n,\Oh_{\PP^n}(k))$), but it can be done after
``rigidification'' -- that is, fixing some discrete structure on $X$. This
is the subject of Invariant Theory in its pre-Hilbert form; however, any
proposed ``geometric'' rigidification depends on the level $k$, and there
is no universal way of doing it. The amazing fact is we can do it in many
cases using the ``classical'' {\em Geometric Quantization Procedure}
(GQP); but for this, we must leave algebraic geometry and go over to
symplectic geometry instead. I would like to call this method the {\em
general theory of theta functions}.

The starting point is that, together with a complex structure $I$ on $X$, a
polarization $L$ gives us a quadruple $(X,\om,L,a_L)$, where $\om$ is the
K\"ahler form and $a_L$ a Hermitian connection on $L$ with curvature form
$F_a=2\pi i\cdot\om$ of Hodge type $(1,1)$, giving the holomorphic
structure on $L$. The pair $(X,\om)$ is a symplectic manifold; we can thus
view it as the {\em phase space} of some classical mechanical system, and
the pair $(L,a_L)$ as {\em prequantization data} of this system.

We should start by recalling the construction of spaces of wave functions
for a pair $(S,\om)$, where $S$ is a smooth symplectic manifold of
dimension $2n$ with a given symplectic form $\om$. To switch on any
quantization procedure, we suppose that the cohomology class $[\om]$ of
the symplectic form is integral, that is, there exists a complex line
bundle $L$ with $c_1(L)=[\om]$. Moreover, suppose that $L$ has a Hermitian
connection $a$ with curvature form $F_a=2\pi i\cdot\om$. Any quadruple of
this type
 \begin{equation}
 (S,\om,L,a)
 \label{eq1.1}
 \end{equation}
is called a {\em prequantization} of the classical mechanical system with
phase space $(S,\om)$.

There are two approaches to the geometric quantization of $(S,\om,L,a)$
(\ref{eq1.1}) (see \cite{A}, \cite{S1} or \cite{W}). We discuss here the
simplest version of these constructions, avoiding questions such as the
choice of metaplectic structures, densities and half densities specifying
geometric conditions on the manifold $S$. (Roughly speaking, $S$ should be
a {\em Calabi--Yau manifold}\/). The usual slogan is that we have to choose
``one half'' of the set of all functions on $S$ using some ``polarization''
conditions. The first approach is as follows:

 \subsubsection*{Complex polarizations}
To define a complex polarization, we give $S$ a complex structure $I$ such
that $S_I=X$ is a K\"ahler manifold with K\"ahler form $\om$. Then the
curvature form of the Hermitian connection $a$ is of type $(1,1)$, hence
for any {\em level} $k\in\Z^+$, the line bundle $L^k$ is a
holomorphic line bundle on $S_I$. Complex quantization provides the space
of {\em wave functions of level $k$}:
 \begin{equation}
 \sH_{L^k}=H^0(S_I,L^k),
 \label{eq1.2}
 \end{equation}
-- that is, the space of {\em holomorphic} sections of $L^k$. Thus a {\em
complex polarization} of $(S,\om,L,a)$ (\ref{eq1.1}) returns to the
algebraic geometry $S=X$ we started from.

In particular, the spaces of wave functions (\ref{eq1.2}) obtained in this
way is the collection of {\em complete linear systems} in the usual sense.
We will suppose $L$ to be an {\em ample} holomorphic line bundle, and in
particular,
 \[
 H^i(S_I,L)=0 \quad\text{for all $i>0$.}
 \]

The second approach is the choice of a real polarization:

\subsubsection*{Real polarizations}
A real polarization of $(S,\om,L,a)$ is a Lagrangian fibration
 \begin{equation}
 \pi\colon S\to B,
 \label{eq1.3}
 \end{equation}
such that
 \[
 \om\rest{\pi^{-1}(b)}=0 \quad\text{for every point $b\in B$,}
 \]
and the fibre $\pi^{-1}(b)$ is a smooth Lagrangian
submanifold for generic $b$.

Thus a mechanical system admits a real polarization if and only if it
is {\em complete integrable}. 

 \begin{rmk}
 Actually, for the ordinary technical tricks of the theory of geometric
quantization to work, we should require that the fibration has regular
geometric behavior (see, for example, \cite{S2}). But beginning with
Guillemin and Sternberg's paper \cite{GS2}, it is reasonable to consider
more general fibrations, namely, {\em real polarizations with
singularities}.
 \end{rmk}

Then restricting $L$ to a Lagrangian fibre gives a flat connection
$a\rest{\text{fibre}}$ or equivalently, a character of the fundamental
group
 \[
 \chi\colon\pi_1(\text{fibre})\to\U(1).
 \]

Let $\sL_{\pi}$ be the sheaf of sections of $L$ that are covariant constant
along fibres. Then we get the space
 \[
 \sH_{\pi}=\bigoplus_{i} H^i(S,\sL_{\pi}).
 \]
In the regular case, \'Sniatycki proved that 
 \[
 H^i(S,\sL_{\pi})=0 \quad\text{for $i \ne n$.}
 \]
 
 \begin{dfn}
 \begin{enumerate}
 \item A fibre of $\pi$ is a {\em Bohr--Sommerfeld} cycle of $(S,\om,L,a)$
if $\chi=1$.
 \item $\BS\subset B$ is the subset of Bohr--Sommerfeld fibres.
 \item $k$-BS${}\subset B$ is the subset of Bohr--Sommerfeld fibres for
$(S,\om,L^k, ka)$.
 \end{enumerate} 
 \end{dfn}

 According to the general theory of real quantizations, we expect to get a
finite number of Bohr--Sommerfeld fibres, and in the regular case,
 \[
 H^n(S,\sL_{\pi})=\bigoplus_{\BS}\C\cdot s_{i},
 \]
where $s_{i}$ is a nonzero covariant constant section of the restriction of
$(L,a)$ to a Bohr--Sommerfeld fibre of the real polarization $\pi$.

In the general case, we can use this to {\em define} the new collection of
spaces of wave functions (of level $k$):
 \begin{equation}
 \sH_{\pi}^k=\bigoplus_{\text{$k$-BS}}\C\cdot s_{i},
 \label{eq1.4}
 \end{equation}
and use special tricks to compare (1.4) with (\ref{eq1.2}).

There is a canonical way of describing the Bohr--Sommerfeld subset. For
this, we must choose special coordinates on $B$, the so-called {\em action
coordinates}, which are part of the {\em action angle} coordinates (see
\cite{A}, \cite{GS1}, \cite{GS2}). Locally around a point $b\in B$, the
action coordinates $c_i$ are given as periods along 1-cycles of the fibre
$\pi^{-1}(b)$ of a 1-form $\al$ such that 
 \begin{equation}
 \dd \al=\om.
 \label{eq1.5}
 \end{equation}
This system of coordinates $\{c_i\}$ is defined up to additive constants
and an {\em integral} linear transformations. Thus, if $B$ is simply
connected, the action coordinates map $B$ locally diffeomorphically to
some open subset 
 \begin{equation}
B_c\subset \R^n_{(c_1, \dots, c_n)}
 \label{eq1.6}
 \end{equation}
with coordinates $\{c_i\}$. If $(0,\dots, 0)$ is a Bohr--Sommerfeld point, then
 \begin{equation}
 \BS=B_c \cap \Z^n
 \label{eq1.7}
 \end{equation}
is the {\em set of integral points} in $B_c$.

Let us return to our collections of spaces of wave functions.
 
 \begin{rmk} An important observation, proved mathematically in a number
of cases, is that the projectivization of the spaces (\ref{eq1.2}) are
given purely by the symplectic prequantization data and do not depend on
the choice of complex structure on $S$. The same is true for the
projectivization of the spaces (\ref{eq1.4}). Moreover, these spaces do
not depend on the real polarization $\pi$ (\ref{eq1.3}), provided that we
extend our prequantization data $(S,\om,L,a,\sF)$ by adding some ``half
density'' $\sF$ (see \cite{GS1}).
 \end{rmk}

Our {\em main problem} is to compare the spaces 
 \[
 \sH_{L^k} \quad \text{and} \quad \sH_{\pi}^k.
 \]
If we are lucky enough to be able to construct a canonical isomorphism
between these spaces, we get a special basis in the space of wave
functions of a complex polarization, and in particular in any ample
complete linear system. To distinguish this basis from others, we call it
the {\em system of theta functions} of level $k$, with ``characteristics''
which are Bohr--Sommerfeld fibres.

Actually, this generalization of the theory of theta functions requires the
final ingredient of the quantization procedure -- the algebra of {\em
observables} represented as an algebra of operators on spaces of wave
functions (like the Heisenberg algebra on spaces of classical theta
functions). We avoid using such algebras in this article, but they
underlie our constructions, so it is reasonable to recall briefly the
general shape of this ingredient.

\subsubsection*{Algebra of observables and its space of states} As a result
of any quantization procedure, we get a $\C^*$-algebra of observables
represented as some algebra $A$ of operators on the spaces of wave
functions (\ref{eq1.2}) or (\ref{eq1.4}). As usual, this algebra is a
noncommutative extension of some commutative $\C^*$-algebra $A_0\subset
A$. For example, if $S=T^*M$ for some mani\-fold $M$ then $A_0$ is the
algebra of continuous complex valued functions, so that $M$ is the {\em
space of maximal ideals} of $A_0$.

A pair $A_0\subset A$ gives us a space $ \sH$ of wave functions
(\ref{eq1.2}) or (\ref{eq1.4}) as the subset of the {\em space of states}.
Recall that a {\em state} is a map:
 \begin{equation}
 \psi\colon A\to\C
 \quad\text{such that} \quad \psi(a^*a)\ge 0 \quad \text{and} \quad \Vert
 \psi \Vert=1.
 \label{eq1.8}
 \end{equation}
The set $\sS(A)$ of all states of $A$ is a convex space and its {\em
boundary elements} are called {\em pure states} (for example, in the
previous example, delta functions of points are pure states). If our
$\C^*$-algebra is represented on $\sH$ by bounded operators then every
vector $\left|\psi\right>$ defines the state as the {\em expectation
value}.

The known strategy to identify spaces (\ref{eq1.2}) and (\ref{eq1.4}) is
to represent both as irreducible representation spaces of some algebra
admitting a {\em unique irreducible representation}. 

The constructions of Berezin, Toeplitz and Rawnsley (see for example
\cite{R}) are extremely useful for our geometric investigations, and we
consider them in \S6.

 \section{Model for our theory: the classical theory of theta functions}

 Let $A$ be a principally polarized Abelian variety of complex dimension
$g$ with flat metric $g$. Then the tangent bundle $TA$ has the standard
constant Hermitian structure (that is, the Euclidean metric, symplectic
form and complex structure $I$). The K\"ahler form $2\pi i\om$ gives a
polarization of degree 1. We fix a {\em smooth} Lagrangian decomposition
of $A$
 \begin{equation}
 A=T^g_+\times T^g_-,
 \label{eq2.1}
 \end{equation}
such that both tori are Lagrangian with respect to $\om$. (In the smooth
category, $A$ is the standard torus $\R^{2g}/\Z^{2g}$ with the standard
constant integral form $\om$, and the decomposition (\ref{eq2.1}) just
consists of putting $\om$ in normal form.) Let $L$ be a holomorphic line
bundle with holomorphic structure given by a Hermitian connection $a$ with
curvature form $F_a=2\pi i\cdot\om$, and $L=\Oh_A(\Theta)$, where $\Theta$
is the classical {\em symmetric} theta divisor. The decomposition
(\ref{eq2.1}) induces a decomposition
 \begin{equation}
 H^1(A,\Z)=\Z^g_+\times\Z^g_-,
 \label{eq2.2}
 \end{equation}
and a Lagrangian decomposition 
 \begin{equation}
 A_k=(T^g_+)_k \times (T^g_-)_k
 \label{eq2.3}
 \end{equation}
of the group of points of order $k$. Any smooth ``irreducible'' $g$-cycle
in $A$ is the image $\fie(T^g)$ of a smooth linear embedding $\fie\colon
T^g\to A$.

 \subsubsection*{Complex quantization}
This is nothing other than the {\em classical theory of theta functions}.
Indeed, the decomposition (\ref{eq2.2}) defines the collection of
compatible {\em theta structures} of every level $k$: the decomposition
(\ref{eq2.3}) defines a Lagrangian decomposition
$A_k=(\Z^g)_k^+\times(\Z^g)_k^-$, and a decomposition of the spaces of
wave functions
 \begin{equation}
 \sH_{L^k}=H^0(A, L^k)=\bigoplus_{w\in (\Z^g)_k^-}\C\cdot\theta_w,
 \quad\text{with}\quad \rk\sH_{L^k}=k^g,
 \label{eq2.4}
 \end{equation}
where $\theta_c$ is the theta function with {\em characteristic} $c$ (see
\cite{Mum}).

The decomposition (\ref{eq2.4}) is given by the following recipe: we
identify the torus $T^g_-$ with the dual torus, and consider vectors
$w\in(T^g_-)_k$ as (periodic) linear differential forms on $T^g_-$.
Applying the symplectic form $\om$ gives a collections of linear vector
fields $\xi_w$ on $A$ parallel to the fibration by the tori $T^g_+$.
Finally, the translations $t_w$ on $A$ obtained as the exponentials of
these vector fields give a finite subgroup of the translations group of
$A$.

Now by choosing $\theta_0\in H^0(A, L^k)$ to be a {\em very symmetric}
section (actually, the section with divisor the sum of all the translates
of the theta divisor $\Theta$ by points of $(T^g_+)_k)$), we get a basis of
$H^0(A, L^k)$:
 \begin{equation}
\{\theta_w=t_w^*(\theta_0)\}.
 \label{eq2.5}
 \end{equation}

 \subsubsection*{Real polarization} The projection of the direct product
(\ref{eq2.1}) gives us a real polarization
 \begin{equation}
\pi\colon A\to T^g_-=B. 
 \label{eq2.6}
 \end{equation}
Remark that in this case the action coordinates (\ref{eq1.6}) are just {\em
flat} coordinates on $T^g_-=B$, and under this identification 
 \begin{equation}
 \kBS=(T^g_-)_k
 \label{eq2.7}
 \end{equation}
is the subgroup of points of order $k$.

 Now we can consider the dual fibration
 \begin{equation}
 \pi'\colon A'=\Pic(A/T^g_-)\to T^g_-=B,
 \label{eq2.8}
 \end{equation}
with fibres
 \[
 (\pi')^{-1} (p)=\Hom(\pi_1(\pi^{-1}(p),\U(1)).
 \]
This fibration admits the section
 \begin{equation}
 s_0\in A \quad\text{with}\quad s_0 \cap (\pi')^{-1}
 (p)=\id\in\Hom(\pi_1(\pi^{-1}(p),\U(1)),
 \label{eq2.9}
 \end{equation}
so that we have a decomposition
 \begin{equation}
 A'=(T^g)'\times T^g_-=B.
 \label{eq2.10}
 \end{equation}

 \begin{rmk} An amazing fact recently proved by Golyshev, Lunts and Orlov
\cite{GLO} is that the $2g$-torus $A'$ is canonically equipped with
 \begin{enumerate}
 \item a symplectic form $\om'$;
 \item a complex structure $I'$.
 \end{enumerate}
 \end{rmk}

Now we can apply geometric quantization to the real polarization
(\ref{eq2.6}) of the phase space $(A,\om, L^k,a_k)$, where $a_k$ is the
Hermitian connection defining the holomorphic structure on $L^k$. Sending
the line bundle $L^k$ to the character of the fundamental group of a fibre
gives a section
 \begin{equation}
s_{L^k}\subset A'=\Pic(A/T^g_-);
 \label{eq2.11}
 \end{equation}
and the Bohr--Sommerfeld subset of $B=T^g_-$ is
 \[
s_0 \cap s_{L^k} \,\subset\, s_0=B=T^g_-\,.
 \]

 Under the identification $s_0=T^g_-=\U(1)^g$, the intersection points
 \[
s_0 \cap s_{L^k}=(\U(1)^g)_k
 \]
are elements of order $k$ in $T^g=\U(1)^g$. We thus get a decomposition
 \begin{equation}
 \sH_{\pi}^k=\bigoplus_{\rho\in\U(1)^g_k}\C\cdot s_{\rho}.
 \label{eq2.12}
 \end{equation}

 \begin{cor} \begin{enumerate} \item $\rk \sH_{L^k}=\rk \sH_{\pi}^k$.
 
 \item Moreover, there exists a canonical isomorphism
 \[
 \sH_{L^k}=\sH_{\pi}^k,
 \]
up to a scaling factor.
 \end{enumerate}
 \end{cor} 

We get already this isomorphism up to the action of the $k^g$-torus
$(\C^*)^{k^g}$ (compare decompositions (\ref{eq2.5}) and (\ref{eq2.12})).
But the canonical iso\-morphism is defined by the action of the Heisenberg
group $H_k$ on holomorphic sections of the line bundle $L^k$ (= the theory
of theta functions, see \cite{Mum}) and the natural extension of the
action of $H_k$ on the collection of Bohr--Sommerfeld orbits. Each of
these representations is irreducible; thus the uniqueness of the
irreducible representation of $H_k$ gives a canonical identification of
these spaces up to scaling.

The functions making up the special bases of these spaces are called {\em
classical theta functions with characteristics of level $k$}.

A real polarization without degenerate fibres such as $\pi$ in
(\ref{eq2.6}) is called {\em regular}. Using more sophisticated techniques
(as in \cite{GS2}) we get a basis of the same type for real polarizations
with degenerate fibres (see Remark after (\ref{eq1.3})). But if we start
with {\em any} polarized K\"ahler manifold $X$, the main question is the
following:
 \begin{quote}
 how to find a real polarization like (\ref{eq1.3}) on $X$ (possibly
 with degenerate fibres)?
 \end{quote}
The amazing fact is that we can do it in many absolutely unpredictable
cases. For example, we now show how to find a real polarization of complex
projective space $\PP^3$. {\bf Warning:} We construct some real
polarization of $\PP^3$, but not a special theta basis in $\PP
H^0(\PP^3,\Oh_{\PP^3}(k))$!

For this, we consider a special presentation of the complex threefold
$\PP^3$ as a real 6-manifold: let $\Si_2$ be a Riemann surface of genus 2.
Then as a 6-manifold,
 \begin{equation}
 \PP^3=\Hom (\pi_1(\Si_2),\SU(2)) /\PU(2)=R_2
 \label{eq2.13}
 \end{equation} 
is the space of classes of $\SU(2)$-representations of the fundamental
group of a Riemann surface of genus 2. 

Thus $\PP^3$ is the first manifold of the collection of manifolds $R_g$. If
we solve the problem of real polarizations of these, we get in particular
a real polarization of $\PP^3$. We do this in the following section, but
we first extend the direct approach by giving a description in terms of
general theories giving rise to these constructions.

 \section{Chern--Simons quantizations of $R_g$}
According to the general procedure, we must present $R_g$ as the classical
phase space of some mechanical system. We begin by recalling the full
steps of this procedure.

A classical field theory on a manifold $M$ has three ingredients:
 \begin{enumerate}
 \item a collection $\sA$ of {\em fields} on $M$, which are geometric
objects such as sections of vector bundles, connections on vector bundles,
maps from $M$ to some auxiliary manifold (the target space) and so on;

 \item an {\em action functional}
 \[
 S \colon \sA\to\C
 \]
which is an integral of a function $L$ (the Lagrangian) of fields;

 \item a collection of observable functionals on the space of fields,
 \[
 \sW \colon \sA\to\C.
 \]
 \end{enumerate}
Our case is the following.

\subsubsection*{Example: Chern--Simons functional} Here $M$ is a
3-manifold, 
 \[
 \sA=\Om^1(M) \tensor \su (2)
 \]
and
 \begin{equation}
S(a)=\frac{1}{8} \pi^2\int_M \tr(a\dd a + \frac{2}{3} a^3). 
 \label{eq3.1}
 \end{equation}
As observable, we can consider a {\em Wilson loop}, given by some
knot $K\subset M$:
 \[
 \sW_K(a)=\tr(\Hol_K(a))
 \]
-- the trace of the holonomy of a connection $a$ around the knot $K$.

Now let 
 \begin{equation}
 R_g=\Hom (\pi_1(\Si_g), \SU(2))/ \PU(2)
 \label{eq3.2}
 \end{equation}
be the space of classes of $\SU(2)$-representations of the fundamental
group of a Riemann surface of genus $g$. This space is stratified by the
subspace of reducible representations
 \begin{equation}
 R_g\red\subset R_g, \quad R_g\irr=R_g-R_g\red.
 \label{eq3.3}
 \end{equation}
To get this space as the phase space of some mechanical system, consider a
compact smooth Riemann surface $\Si$ of genus $g > 1$ and the trivial
Hermitian vector bundle $E_h$ of rank 2 on it. As usual, let $\sA_h$ be the
affine space (over the vector space $\Om^1(\End E_h)$) of Hermitian
connections and $\sG_h$ the Hermitian gauge group. This space admits a
stratification:
 \[
 \sA_h\red\subset \sA_h
 \]
where the left-hand side is the subset of reducible connections. As usual,
let
 \[
\sA_h\irr=\sA_h-\sA_h\red.
 \]

Sending a connection to its curvature tensor defines a $\sG_h$-equivariant
map
 \begin{equation}
F \colon \sA(E_h)\to\Om^2(\End E_h)=\Lie(\sG_h)^*
 \label{eq3.4}
 \end{equation}
to the coalgebra Lie of the gauge group. 

We can consider this map as the moment map with respect to the action of
$\sG_h$. The subset
 \begin{equation}
F^{-1}(0)=\sA_F
 \label{eq3.5}
 \end{equation}
is the subset of flat connections and 
 \[
\sA_F\irr=\sA_F \cap \sA_h\irr
 \]
the subspace of irreducible flat connections.

For a connection $a\in \sA_F$ and a tangent vector to $\sA_h$ at $a$
 \[
 \om\in\Om^1(\End E_h)=T \sA_h,
 \]
we have
 \begin{equation}\om\in T \sA_F \iff \na_a (\om)=0.
 \label{eq3.6}
 \end{equation} 
The trivial vector bundle $E_h$ admits the trivial connection $\theta$,
which is interesting and important from many points of view, and it
provides in particular the possibility of identifying $\sA_h$ with
$\Om^1(\End E_h)$ by sending a connection $a$ to the form $a-\theta$. We
will identify forms and connections in this way.

 The space $\sA_h=\Om^1(\End E_h)$ is the {\em collection of fields} of
YM-QFT with the {\em Yang--Mills functional}
 \[
 S(a)=\int_{\Si} |F_a|^2.
 \]
Thus $\sA_F$ is a {\em classical phase space}, that is, the space of
solutions of an {\em Euler--Lagrange equation} $\de S(a)=0$.

There exists a symplectic structure on the affine space $\sA_h$, induced by
the canonical 2-form given on the tangent space $\Om^1(\End E_h)$ at a
connection $a$ by the formula 
 \begin{equation}
 \Om_0 (\om_1,\om_2)=\int_{\Si} \tr (\om_1 \wedge\om_2). 
 \label{eq3.7}
 \end{equation} 
This form is $\sG_h$-invariant, and its restriction to $\sA_F\irr$ is
degenerate along $\sG_h$-orbits: at a connection $a$, for a tangent vector
$\om\in\Om^1(\End E_h)$, we have
 \[
 \om\in T\sG_h \iff \om=\na_a \fie \quad\text{for $
 \fie\in\Om^0(\End E_h)=\Lie(\sG_h)^*$,}
 \]
and
 \[
 \int_{\Si} \tr (\na_a \fie \wedge\om)=\int_{\Si}\tr(\fie\wedge\na_a\om)=0.
 \]
Hence
 \begin{equation}
 \om\in T \sA_F \iff \Om_0 (\na_a \fie,\om)=0.
 \label{eq3.8}
 \end{equation}

Interpreting (\ref{eq3.4}) as a moment map and using symplectic reduction
arguments, we get a nondegenerate closed symplectic form $\Om$ on the space
 \[
 \sA_F /\sG_h=R_g
 \]
of classes of $\SU(2)$-representations of the fundamental group of the
Riemann surface, and a stratification of this space. The form $\Om$
defines a symplectic structure on $R_g\irr$ and a symplectic orbifold
structure on $R_g$.

On the other hand, the form $\Om_0$ on $\sA_h$ is the differential of the
1-form $D$ given by the formula
 \begin{equation}
D(\om)=\int_{\Si} \tr((a) \wedge\om).
 \label{eq3.9}
 \end{equation}
We consider this form as a unitary connection $A_0$ on the trivial
principal \hbox{$\U(1)$-bundle} $L_0$ on $\sA_h$.

To descend this Hermitian bundle and its connection to the orbit space,
one defines the $\Theta$-cocycle (or $\Theta$-torsor) on the trivial line
bundle (see \cite{RSW}). This cocycle is the $\U(1)$-valued function
$\Theta$ on $\sA_h\times\sG_h$ defined as follow: for any triple
$(\Si,a,g)$ where $(a,g)\in\sA_h\times\sG_h$, we can find a triple
$(Y,A,G)$ where $Y$ is a smooth compact 3-manifold, $A$ a
$\SU(2)$-connection on the trivial vector bundle $\sE$ on $Y$ and $G$ a
gauge transformation of it, such that
 \[
 \p Y=\Si, \quad a=A \rest{\Si} \quad \text{and} \quad g=G\rest{\Si}.
 \]
Then
 \begin{equation}
 \Theta (a, g)=e^{i(\CS(A^G)-\CS(A))}.
 \label{eq3.10}
 \end{equation}

Recall that the Chern--Simons functional on the space $\sA (\sE_h)$ of
unitary connections on the trivial vector bundle is given by the formula
 \begin{equation}
 \CS_Y (a_0 +\om)=\int_Y \tr\left(\om \wedge F_{a_0}-\frac{2}{3}
 \om\wedge\om\wedge\om\right).
 \label{eq3.11}
 \end{equation} 
It can be checked that the function (\ref{eq3.10}) does not depend on the
choice of the triple $(Y,A,G)$ (see \cite{RSW}, \S2). 

The differential of $\Theta$ at $(a, g)$ is given by the formula
 \begin{multline}
 \dd\Theta(\om,\fie)= \\
 \frac{\pi i}{4}\,\Theta\int_{\Si}\Bigl(\tr(g^{-1}\dd g\wedge g^{-1}\om
 g)-\tr(a\wedge\na_{a^g}\fie)+2\tr(F_{a^g} \wedge \fie) \Bigr),
 \label{eq3.12}
 \end{multline}
where $\om\in\Om^1(\End E_h)$ and $\fie\in\Om^0(\End E_h)=\Lie(\sG_h)$.

But the restriction of this differential to the subspace of flat
connections is much simpler:
 \begin{equation}
 \dd \Theta (\om, \fie)=\frac{\pi i}{4}\,\Theta\int_{\Si}\tr(g^{-1}\dd g
 \wedge g^{-1}\om g),
 \label{eq3.13}
 \end{equation} 
and is independent of the second coordinate.

That this function is in fact a cocycle results from the functional
equation
 \begin{equation}
\Theta (a, g_1 g_2)=\Theta (a, g_1) \Theta (a^{g_1}, g_2). 
 \label{eq3.14}
 \end{equation}
Using this function as a torsor $\sA_h \times_{\Theta}\U(1)$ we get a
principal $\U(1)$-bundle $S^1(L)$ on the orbit space $\sA_h/\sG_h$:
 \begin{equation}
S^1(L)=(\sA_h\times S^1)/\sG_h,
 \label{eq3.15}
 \end{equation}
where the gauge group $\sG_h$ acts by
 \[
g(a, z)=(a^g, \Theta(a,g) z),
 \]
or the line bundle $L$ with a Hermitian structure. 

Following \cite{RSW}, let us restrict this bundle to the subspace of flat
connections $\sA_F$. Then one can check that the restriction of the form
$D$ (\ref{eq3.9}) to $\sA_F$ defines a $\U(1)$-connection $A_{\CS}$ on the
line bundle $L$.

By definition, the curvature form of this connection is
 \begin{equation}
F_{A_{\CS}}=i\cdot\Om.
 \label{eq3.16}
 \end{equation}

Thus the quadruple
 \begin{equation}
(R_g,\Om, L, A_{\CS})
 \label{eq3.17}
 \end{equation}
is a {\em prequantum system} and we are ready to switch on the Geometric Quantization Procedure.

\section{Complex polarization of $R_g$} 

The standard way of getting a complex polarization is to give a Riemann
surface $\Si$ of genus $g$ a conformal structure $I$. We get a complex
structure on the space of classes of representations $R_g$ as follows: let
$E$ be our complex vector bundle and $\sA$ the space of all connections on
it. Every connection $a\in\sA$ is given by a covariant derivative
$\na_a\colon\Ga(E)\to\Ga(E\tensor T^*X)$, a first order differential
operator with the ordinary derivative $\dd$ as the principal symbol and a
complex structure gives the decomposition $\dd=\p+\dbar$, so any covariant
derivative can be decomposed as $\na_a=\p_a+\dbar_a$, where
$\p_a\colon\Ga(E)\to\Ga(E\tensor\Om^{1,0})$ and
$\dbar_a\colon\Ga(E)\to\Ga(E\tensor\Om^{0,1})$. Thus the space of
connections admits a decomposition 
 \begin{equation}
 \sA=\sA'\times \sA'',
 \label{eq4.1}
 \end{equation}
where $\sA'$ is an affine space over $\Om^{1,0}(\End E)$ and $\sA''$ an
affine space over $\Om^{0,1}(\End E)$.

The group $\sG$ of all automorphisms of $E$ acts as the group of gauge
transformations, and the projection $\pr\colon \sA\to\sA''$ to the space
$\sA''$ of $\dbar$-operators on $E$ is equivariant with respect to the
$\sG$-action.

Giving $E$ a Hermitian structure $h$, we get the subspace $\sA_h\subset
\sA$ of Hermitian connections, and the restriction of the projection $\pr$
to $\sA_h$ is one-to-one. Under this Hermitian metric $h$, every element
$g\in\sG$ gives an element $\gbar=(g^*)^{-1}$ such that
 \[
 \gbar=g \iff g\in\sG_h.
 \]
Now for $g\in\sG$, the action of $\sG$ on the component $\sA''$ is
standard:
 \[
 \dbar_{g(a)}=g\cdot \dbar_a\cdot g^{-1}=\dbar_a-(\dbar_a g)\cdot g^{-1};
 \]
and the action on the first component $\sA'$ of $\p$-operators is
 \[
 \p_{g(a)}=\gbar\cdot\p_a\cdot\gbar^{-1}=\p_a-((\dbar_a g)\cdot g^{-1})^*.
 \]
It is easy to see directly that the action just described preserves unitary
connections:
 \begin{equation}
 \sG(\sA_h)=\sA_h,
 \label{eq4.2}
 \end{equation}
and that the identification $\sA_h=\sA$ is equivariant with respect to
this action.

It is easy to see that $\dbar_a^2\in\Om^{0,2}(\End E)=0$. Thus the orbit
space
 \begin{equation}
\sA'' /\sG=\bigcup \sM_i 
 \label{eq4.3}
 \end{equation} 
is the union of all components of the moduli space of topologically
trivial \hbox{$I$-holomorphic} bundles on $\Si_I$. (This union doesn't
admit any good structure, as it contains all unstable vector bundles).
Finally, the image of $\sA_F\in \sA_h$ is the component $\sM\semis$ of
maximal dimension ($3g-3$) of s-classes of semistable vector bundles. Thus
by classical technique of GIT of Kempf--Ness type we get:

 \begin{prop} {\em (Narasimhan--Seshadri)} 
 \[
 R_{\Si}=R_g=\sM\semis.
 \] 
 \end{prop}

 \begin{prop} The form $F_{A_{\CS}}$ (\ref{eq3.16}) is a $(1,1)$-form and 
the line bundle $L$ admits a unique holomorphic structure compatible with
the Hermitian connection $A_{\CS}$.
 \end{prop}

On the other hand, a complex structure $I$ on $\Si$ defines a K\"ahler
metric on $\sM\semis$ (the so-called Weyl--Petersson metric) with K\"ahler
form
 \begin{equation}\om_{\mathrm{WP}}=i F_{A_{\CS}}=i\cdot\Om.
 \label{eq4.4}
 \end{equation}
This metric defines the Levi-Civita connection on the complex tangent
bundle $T \sM\semis$, and hence a Hermitian connection $A_{\LC}$ on the
line bundle
 \begin{equation}
\det T \sM\semis=L^{\tensor 4},
 \label{eq4.5}
 \end{equation} 
and a Hermitian connection $\frac{1}{4}A_{\LC}$ on $L$ compatible with the
holomorphic structure on $L$. Thus we have

 \begin{prop} \label{prop4.3}
 \[ 
 \frac{1}{4} A_{\LC}=A_{\CS}.
 \] 
 \end{prop}

Finally, considering $\sM\semis$ as a family of $\dbar$-operators, we get
the Quillen determinant line bundle $L$ having a Hermitian connection
$A_Q$ with curvature form
 \begin{equation}
F_{A_Q}=i\cdot\Om.
 \label{eq4.6}
 \end{equation} 
Hence we can extend the equality of Proposition~\ref{prop4.3}:
 \begin{prop} 
 \[ 
 \frac{1}{4} A_{\LC}=A_{\CS}=A_Q.
 \] 
 \end{prop}

Summarizing, the result of the complex quantization procedure of the
prequantum system (\ref{eq3.17}) can be considered to be the spaces of
wave functions of level $k$, that is, the spaces of $I$-holomorphic
sections
 \begin{equation} 
 \sH_{L^k}=H^0 (L^k)
 \label{eq4.7}
 \end{equation} 

One knows that this system of spaces and monomorphisms is related to the
system of representations of $\fsl(2,\C)$ in the
Weiss--Zumino--Novikov--Witten model of CQFT. Namely, for a half integer
$i$, consider the irreducible representation $V_i$ of dimension $2i+1$ of
$\fsl(2,\C)$. The tensor product of two such representations is given by
the Clebsch--Gordan rule
 \begin{equation}
 V_i\tensor V_j=V_{i+j}\oplus V_{i+j-1}\oplus\dots\oplus V_{i-j}\quad
 \text{for $i\ge j$,}
 \label{eq4.8}
 \end{equation}
and the level of $V_i$ is $2i$. 

Then the {\em fusion
ring}\linebreak[3] $R_k(\fsl(2,\C))$ of level $k$ is the quotient
 \begin{equation}
R_k(\fsl(2,\C))=R(\fsl(2,\C))/\Span{V_{(k+1)/2}}
 \label{eq4.9}
 \end{equation}
of the {\em representation ring} $R(\fsl(2,\C))$ by the ideal generated by
$V_{(k+1)/2}$.

Moreover, every character of the ring $R(\fsl(2,\C))$ is given by a
complex number $z\in\C$ which we can consider as a diagonal $2\times2$
matrix $\diag(iz,-iz)$. This matrix acts on $\fsl(2,\C)$ and $V_i$ and 
 \begin{equation}
\chi_z(V_i)=\tr(\exp(\diag(iz,-iz)))=\frac{\sin((2i+1)z)}{\sin z}.
 \label{eq4.10}
 \end{equation} 
Thus 
 \begin{equation}
\chi_z(V_i)=0 \iff z=\frac{n \pi}{k+2}\quad \text{for $1 \le n \le 2i+1$.}
 \label{eq4.11}
 \end{equation}
In these terms we get:
 \begin{equation}
 \sH_{L^k}=\frac{(k+2)^{g-1}}{2^{g-1}} \sum_{n=1}^{k+1}
 \frac{1}{(\sin(\frac{n\pi}{k+2}))^{2g-2}}\,.
 \label{eq4.12}
 \end{equation}

See \cite{B} for a mathematical derivation of this formula.

\section{Real polarization of $R_g$}

 The collection of real polarizations of the prequantum system 
 \[
 (R_g,\Om, L, A_{\CS}) 
 \]
is given in a very geometric way in the set-up of perturbation theory of
\hbox{3-dimensional} Chern--Simons theory. The crucial point is a {\em
trinion decomposition} of a Riemann surfaces, given by a choice of a
maximal collection of disjoint, noncontractible, pairwise nonisotopic
smooth circles on $\Si$. An isotopy class of such a collection of circles
is called a {\em marking} of the Riemann surface. It is easy to see
(\cite{HT}) that any such system contains $3g-3$ simple closed circles 
 \begin{equation} 
C_1, \dots, C_{3g-3}\subset \Si_g,
 \label{eq5.1}
 \end{equation} 
and the complement is the union
 \begin{equation}
 \Si_g-\{C_1, \dots, C_{3g-3}\}=\bigcup_{i=1}^{2g-2} P_i
 \label{eq5.2}
 \end{equation}
of $2g-2$ trinions $P_i$, where every trinion is a 2-sphere with 3
disjoint discs deleted:
 \[
 P_i=S^2 \setminus \bigl(D_1\cup D_2\cup D_3\bigr) \quad \text{with}
 \quad\Dbar_i\cap\Dbar_j=\emptyset \quad \text{for} \quad i \ne j.
 \]
A collection $\{C_i\}$ with these conditions is called a {\em trinion
decomposition} of $\Si$. The invariant of such a decomposition by marking
class is given by its {\em $3$-valent dual graph} $\Ga(\{C_i\})$,
associating a vertex to each trinion $P_i$, and an edge linking $P_i$ and
$P_j$ to a circle $C_l$ (\ref{eq5.1}) such that 
 \[
 C_l\,\subset\,\p P_i \cap\p P_j.
 \]

If we fix the isotopy class of a trinion decomposition $\{C_i\}$, we get a
map
 \begin{equation}
\pi_{\{C_i\}} \colon R_g\to \R^{3g-3}
 \label{eq5.3}
 \end{equation} 
with fixed coordinates $(c_1, \dots, c_{3g-3})$ such that 
 \[
 c_i (\pi_{\{C_i\}} (\rho))=\frac{1}{\pi}\,\cos^{-1}\bigl(\frac{1}{2}\tr
\rho([C_i])\bigr)\in [0, 1].
 \]

We see that 

 \begin{prop} 
 \begin{enumerate}
 \item The map $\pi_{\{C_i\}}$ is a real polarization of the system
$(R_g,k\cdot\om,L^k,k\cdot A_{\CS})$.

 \item The coordinates $c_i$ are {\em action coordinates} for this
Hamiltonian system (see (\ref{eq1.6}) and \cite{D}).
 \end{enumerate}
 \end{prop}

These functions $c_i$ are continuous on all $R_g$ and smooth over $(0,1)$.
Moreover, Goldman \cite{G} constructed $\U(1)$-actions on the open dense
set
 \[
 U_i=c_i^{-1}(0,1) 
 \]
for which the function $c_i$ is the {\em moment map}, and all these
$\U(1)$-actions commute with each other. Hence we get:

 \begin{prop}
 \begin{enumerate}
 \item The restriction 
 \[
 \pi_{\{C_i\}}\rest{\bigcap_i U_i}\colon\bigcap_i U_i\to (0, 1)^{3g-3}
 \]
is the moment map for the $\U(1)^{3g-3}$-action on $\bigcap_i U_i$.

 \item The image of $R_g$ under $\pi_{\{C_i\}}$ is a convex polyhedron 
 \begin{equation}
 I_{\{C_i\}}\subset [0, 1]^{3g-3}.
 \label{eq5.4}
 \end{equation} 

 \item The symplectic volume of $R_g$ equals the Euclidean volume of
$I_{\{C_i\}} $:
 \begin{equation}\int_{R_g}\om^{3g-3}=\Vol I_{\{C_i\}}.
 \label{eq5.5}
 \end{equation}

 \item The expected number of Bohr--Sommerfeld orbits of the real
polarization $\{C_i\}$
 \begin{equation}
 N_{\BS}(\pi_{\{C_i\}},R_g,\om,L,A_{\CS}) 
 \label{eq5.6}
 \end{equation}
equals the number of half integer points in the polyhedron $I_{\{C_i\}}$,
and 
 \begin{equation}
 \lim_{k\to\infty} k^{3-3g}\cdot N_{\kBS}=\int_{R_g}\om^{3g-3}
 =\Vol I_{\{C_i\}}.
 \label{eq5.7}
 \end{equation}
 \end{enumerate}

 \end{prop}

{From} the combinatorial point of view, the number $N_{\BS}$, or more
generally the numbers $N_{\kBS}$ of $k$-BS fibres, is determined as
follows: consider functions
 \begin{equation}
 w\colon\{C_i\}\to\frac{1}{2k}\{0,1,2,\dots,k\}
 \label{eq5.8}
 \end{equation}
on the collection of edges of the 3-valent graph $\Ga (\{C_i\})$ to the
collection of $\frac{1}{2k}$ integers, such that, for any three edges $C_l,
C_m, C_n$ meeting at a vertex $P_i$, the following 3 conditions hold:
 \begin{enumerate}
 \item $w(C_l) + w(C_m) + w(C_n)\in \frac{1}{k}\cdot \Z$;

 \item $w(C_l) + w(C_m) + w(C_n) \le 1$;

 \item for any ordering of the triple $C_l, C_m, C_n$,
 \begin{equation}
 |w(C_l)-w(C_m)|\le w (C_n) \le w(C_l) + w(C_m).
 \label{eq5.9}
 \end{equation}
 \end{enumerate}
Such a function $w$ is called an {\em admissible integer weight of level}
$k$ on the graph $\Ga(\{C_i\})$.

 \begin{prop}
 \begin{enumerate}
 \item The number $|W_g^k|$ of admissible weights of level $k$
is independent of the graph $\Ga (\{C_i\})$;

 \item
 \begin{equation}
 |W_g^k|=N_{\kBS}.
 \label{eq5.10}
 \end{equation}
 \end{enumerate}
 \end{prop} 

The conditions (\ref{eq5.9}) are called {\em Clebsch--Gordan conditions}
for $\fsl(2,\C)$, for obvious reasons. We can view the space of all real
functions with values in $[0,1]$ subject to these conditions to get a
complex $I_{\{C_i\}}$.

 \begin{rmk} Following this combinatorial approach, we can construct a two
dimensional complex $Y_g$: the set of vertices is the set of all dual
graphs associated with all types of markings of $\Si$. Two vertices are
joined by an edge if and only if the two graphs are related by an {\em
elementary fusion operation}. The $2$-cells correspond to {\em pentagons},
and so on (see \cite{MS}). The topology of this complex reflects the
combinatorial properties of real polarizations of this type.

The geometric meaning of this combinatorial description is as follows:
consider the space $\R^{3g-3}$ with action coordinates $c_i$
(\ref{eq5.3}). This space contains the integer sublattice
$\Z^{3g-3}\subset \R^{3g-3}$, and we can consider the {\em action torus}:
 \begin{equation}
 T^A=\R^{3g-3} / \Z^{3g-3}.
 \label{eq5.11}
 \end{equation}

In particular, we get a map
 \begin{equation}
 \pi_A \colon R_g\to T^A
 \label{eq5.12}
 \end{equation}
which glues at most points of the boundary of $I_{\{C_i\}}$.
 \end{rmk} 

Now every integer weight $w$ (\ref{eq5.8}) satisfying (1) and (2), but a
priori without the Clebsch--Gordan conditions (\ref{eq5.9}), defines a
point of order $2k$ on the action torus
 \[
 w\in T^A_{2k}.
 \]
In particular, the collection $W_g^k$ of admissible integer weights
(subject to (\ref{eq5.9})) can be considered as a subset of points of
order $2k$ on the action torus:
 \begin{equation}
 W_g^k\subset T^A_{2k}.
 \label{eq5.13}
 \end{equation}
On the other hand, every vector $w\in \R^{3g-3}$ can be interpreted as a
{\em differential $1$-form} on $\R^{3g-3}$, and by the usual construction
using the symplectic form $\Om$, this defines a vector field $\xi_w$
tangent to the fibres of $\pi$. Integrating such vector fields defines the
collection of transformations
 \begin{equation}
 \{t_w\}=e^{\xi_w}\subset\Diff^+ (R_g).
 \label{eq5.14}
 \end{equation}

These transformations preserve the curvature form $A_{\CS}$ of the
connection. Thus (because $R_g$ is simple connected), there exists a
collection of gauge transformations $\al_w\in\sG_L$ of $L$ such that
 \begin{equation}
 (t_w)^*(A_{\CS})=A_{\CS}^{\al_w}.
 \label{eq5.15}
 \end{equation}
We can view such gauge transformations as $\U(1)$-{\em torsors}, just as
in describing the formulas for classical theta functions for Abelian
varieties in \S2.

Moreover, if $R_{\Si}$ is given the K\"ahler structure induced from $\Si$
and $s\in H^0(R_{\Si}, L^k)$ is a holomorphic section, then we have the
following.
 \begin{prop}
 \begin{equation}
(t_w)^*(s)\in H^0(R_{\Si}, L^k) 
 \label{eq5.16}
 \end{equation}
is also a holomorphic section.
 \end{prop}
 \begin{cor} If $s_0$ is a {\em sufficiently symmetric} holomorphic section
of $L^k$, then the system
 \begin{equation}
\{s_w=t_w^*(s_0)\}\subset H^0(R_{\Si}, L^k) 
 \label{eq5.17}
 \end{equation}
is a special theta basis of\/ {\em some subspace} of $H^0(L^k)$.
 \end{cor}
Comparing (\ref{eq5.17}) and (\ref{eq2.6}), we see that the recipe to
construct the theta basis is the same as for Abelian varieties with the
action space $T^A$ (see \S2) but instead of the full collection $T^A_{2k}$
of points of order $2k$, we only use the subset $W_g^k\subset T^A_{2k}$.

In our realistic situation, the prequantum system $(R_g,\Om,L,A_{\CS})$ is
far from the regular ``theoretical'' case. But in the fundamental papers
\cite{JW1} and \cite{JW2} there is a well-defined correction to the
``theoretical'' situation. Here we only explain what we must do at a
maximally degenerate Bohr--Sommerfeld fibre. We get proofs of the central
statements of Proposition~5.4 and Corollary~5.1 by a quite fruitful
method: we give new definitions making the statements almost obvious. We
do this in the following special section.

\subsubsection*{Unitary Schottky representations}

Every oriented trinion $P_i$ defines a {\em handle\/} $\HP_i$, and all
these handles glue together to give a {\em handlebody\/} $H_{\{C_i\}}$, a
compact 3-manifold such that
 \begin{equation}\p H_{\{C_i\}}=\Si.
 \label{eq5.18}
 \end{equation}
We get a surjection
 \begin{equation}
\fie \colon \pi_1(\Si)\to \pi_1(H_{\{C_i\}}),
 \label{eq5.19}
 \end{equation}
which defines the subspace 
 \begin{equation}
B_{\{C_i\}}=\bigl\{\rho\in R_g \bigm| \rho \rest{\ker \fie}=1\bigr\}.
 \label{eq5.20}
 \end{equation}

 \begin{prop}
 \begin{enumerate}
 \item $B_{\{C_i\}}$ is a Lagrangian subspace of $R_g$.

 \item More precisely, it is a fibre of the real polarization
$\pi_{\{C_i\}}$ (\ref{eq5.3}):
 \begin{equation}
B_{\{C_i\}}=\pi_{\{C_i\}}^{-1} (1,\dots,1).
 \label{eq5.21}
 \end{equation}

 \item Moreover, $B_{\{C_i\}}$ is a Bohr--Sommerfeld orbit of
$\pi_{\{C_i\}}$.
 \end{enumerate}
 \end{prop}

This Lagrangian subspace $B_{\{C_i\}}$ is singular:
 \[
 B_{\{C_i\}}\red=B_{\{C_i\}} \cap R_g\red=\Sing B_{\{C_i\}};
 \]
and $B_{\{C_i\}}\irr=B_{\{C_i\}} \cap R_g\irr$ is a nonsingular Lagrangian
subvariety.

Under the identification of $R_g$ with the moduli space of s-classes of
semi\-stable vector bundles on the algebraic curve $\Si$, the subspace
$B_{\{C_i\}}$ is called the subset of {\em unitary Schottky subbundles}.

Obviously for this Bohr--Sommerfeld fibre $w_{\US}\in T^A_k$, we must use
a special description of the symplectomorphism $t_{w_{\US}}$. This was
done in the papers \cite{JW1} and \cite{JW2}. 

Returning to the general geometric quantization procedure and summarizing
these results, we get two spaces of wave functions: complex quantization
gives the spaces
 \[
 \sH_{\Si}^k=H^0(L^k)
 \]
of $I$-holomorphic sections of $L^k$, and real quantization gives the
direct sum
 \begin{equation}
 \sH_{\pi}^k\,=\sum_{k\text{-BS fibres}}\!\C\cdot s_{\kBS}
 \label{eq5.22}
 \end{equation}
of lines generated by covariant constant sections of restrictions of our
prequantum line bundle $(L^k,k\cdot A_{\CS})$ to the Bohr--Sommerfeld
fibres of $\pi$ of level $k$.
 
The amazing fact is the following:

 \begin{prop} For any level $k$, any complex Riemann surface $\Si$, and any
trinion decomposition $\{C_i\}$ with the real polarization $\pi$ of $R_g$
we have
 \begin{equation}
\rk \sH_{\Si}=H^0(L^k)=\rk \sH_{\pi},
 \label{eq5.23}
 \end{equation}
and these ranks can be computed by the Verlinde calculus.
 \end{prop}

 \begin{cor} Our construction gives a distinguished theta basis of the
first space $H^0(L^k)$.
 \end{cor}

This isomorphism between spaces of wave functions underlies all the
``modular'' behavior of gauge theory invariants in dimensions 2, 3 and 4.

The final ``classical'' question concerns the Fourier decomposition of our
non-Abelian theta functions $s_w$ (5.17). It can be done using the Fourier
decomposition along coordinates $\{c_i\}$ of the action torus $T^A$
(\ref{eq5.11}) twisting by the system of torsors $\{\al_w\}$ (5.15).
Roughly speaking, the theta functions $s_w$ (5.17) are {\em truncated}
theta functions on the $(6g-6)$-dimensional ``Fourier torus''
 \[
 T_F=\U(1)^{3g-3} \times T^A.
 \]
Namely all coefficients of Fourier decompositions not satisfying the
Clebsch--Gordan conditions (\ref{eq5.9}) must go to zero. Can this
condition be interpreted in terms of the heat equation?

\section{Other definition of a theta basis}

We must first recall the main constructions of GQP. Let $h$ be the
Hermitian form on $L$, and
 \begin{equation}
\mu=\frac{1}{(3g-3)!}\om^{3g-3}
 \label{eq6.1}
 \end{equation}
the volume form on $R_g$. Then we have a scalar product and norm on the
space $\Ga (L^k)$ of global differentiable sections of $L^k$:
 \begin{equation}
 \Span{s_1,s_2}=\int_{R_g} h(s_1, s_2)\cdot \mu \quad \text{and} \quad
 \Vert s \Vert=\sqrt{\Span{s,s}}\,.
 \label{eq6.2}
 \end{equation}
Let $L^2(L^k)$ be the $L^2$-completion of $\Ga(L^k)$ and
 \begin{equation}
P_k \colon L^2(L^k)\to H^0(L^k)
 \label{eq6.3}
 \end{equation}
the orthogonal projection to the finite dimensional subspace of {\em
holomorphic} sections $H^0(L^k)\subset L^2(L^k)$.

The ring $C^{\infty}(R_g)$ of smooth functions on $R_g$ acts on $L^2(L^k)$
by multi\-plication $s\to f\cdot s$, and acts on the space $H^0(L)$ as a
{\em Toeplitz operator}:
 \begin{equation}
T_f=P \odot f\in\End(H^0(L^k));
 \label{eq6.4}
 \end{equation}
the map $C^{\infty}(R_g)\to\End(H^0(L^k))$ is called the {\em
Berezin--Toeplitz map}.

Now, let 
 \begin{equation}
p \colon L^*\to R_g
 \label{eq6.5}
 \end{equation}
 be the principal $\C^*$-bundle of $L$. Every point $x\in L^k$ defines a
linear form
 \begin{equation}
l_x \colon H^0(L^k)\to\C, \quad\text{given by}\quad s(p(x))=l_x (s)\cdot x,
 \label{eq6.6}
 \end{equation}
and the {\em coherent state} vector $s_x\in H^0(L^k)$ associated to $x$,
which is uniquely determined by the equation
 \begin{equation}
\Span{s_x,s}=l_x(s).
 \label{eq6.7}
 \end{equation}
Thus we get a map
 \begin{equation}
 \fie_k \colon R_g\to \PP H^0(L^k),
 \label{eq6.8}
 \end{equation}
which is nothing other than the Hermitian conjugate of the standard
algebraic geometric map by a complete linear system to the dual space,
because of the equality 
 \[
 s_{\al\cdot x}=\albar^{-1}\cdot s_x \quad \text{for $\al\in\C^*$.}
 \] 
Now, following John Rawnsley, we can define {\em coherent projectors}
$P_{p(x)}$ and the {\em Rawnsley epsilon function} $\ep\colon R_g\to\R^+$
in such a way that:
 \begin{equation}
 \ep(p(x))=|x|^2\cdot \Span{s_x,s_x} \quad \text{and} \quad
 h(s_1, s_2)_{p(x)}=\ep (p(x))\cdot \Span{s_1, P_{p(x)} s_2}. 
 \label{eq6.9}
 \end{equation}
Since $\ep > 0$ we can modify the old measure on $R_g$:
 \begin{equation}
\mu_{\ep}=\ep\cdot \mu,
 \label{eq6.10}
 \end{equation}
where $\mu$ is (6.1). This measure gives an integral representation of
Toeplitz operators:
 \begin{equation}
T_f(s)=\int_{R_g} f(p(x))\cdot P_{p(x)}(s)\cdot \mu_{\ep}.
 \label{eq6.11}
 \end{equation}

Up to now, we have been working with a {\em complex polarization}. Let us
return to the real polarization $\pi$ (\ref{eq5.3}). We can identify the
target real space $\R^{3g-3}$ of $\pi$ with the dual space 
 \[
\R^{3g-3}=(\R^{3g-3})^*
 \]
and we can consider our vectors $w\in W^k_g\subset T^A_{2k}$ as {\em
linear functions} on the target space $\R^{3g-3}$. Thus we get a
collection of functions
 \begin{equation}
\pi^* w \colon R_g\to \R
 \label{eq6.12}
 \end{equation}
and a collection of Toeplitz operators
 \begin{equation}
T_{\pi^* w}\in\End(H^0(L^k)).
 \label{eq6.13}
 \end{equation}
Let us choose one (very symmetric) section $s_0$ in the following way: for
$k=1$, the space $H^0(L)$ is the space of ordinary theta functions (see,
for example, \cite{BL}) and every semistable bundle $E$ defines a theta
divisor
 \[
 \Theta_E=\bigl\{L\in\Pic_{g-1}(\Si)\bigm| h^0(E\tensor L)>0\bigr\}.
 \] 
Let $s_E$ be the section with this divisor as its zero set. Then one has
the section
 \begin{equation}
 s_0=s_{\Oh_{\Si}\oplus\Oh_{\Si}}^k\in H^0(R_{\Si},L^k),
 \label{eq6.14}
 \end{equation}
and the collection of sections
 \begin{equation}
 \{T_{\pi^* w}(s_0)=s_w\}\subset H^0(R_{\Si},L^k). 
 \label{eq6.15}
 \end{equation}

Using the integral representation (\ref{eq6.11}), Rawnsley's localization,
and Proposition 5.6, we get immediately
 \begin{thm} The sections $\{T_{\pi^* w}(s_0)=s_w\}$ form a basis of
$H^0(R_\Si, L^k)$.
 \end{thm}

The reader not wishing to check the following statement may take
(\ref{eq6.15}) as the {\em definition} of the theta basis:

 \begin{prop} The basis (5.17) coincides with the basis (\ref{eq6.15}).
 \end{prop}

Thus, we do indeed get a generalization of theta functions.

\section{What next?} The theory of theta functions outlined above is just
a small sample of the applications of the Geometric Quantization Procedure
to algebraic geometry. Here we extend the list, mentioning applications
which are natural generalizations of the above constructions.

\subsubsection*{Generalization to vector bundles of higher rank}
This construction is new even in the classical set-up. Let us return to a
real polarization of an Abelian variety $\pi\colon A\to T^g_-=B$, and its
dual fibration
 \[
 \pi'\colon A'=\Pic(A/T^g_-)\to T^g_-=B,
 \]
with fibres
 \[
 (\pi')^{-1} (p)=\Hom(\pi_1(\pi^{-1}(p)),\U(1))
 \]
and section
 \[
 s_0\in A' \quad\text{with}\quad s_0 \cap (\pi')^{-1}
 (p)=\id\in\Hom(\pi_1(\pi^{-1}(p)),\U(1)).
 \]
Every stable holomorphic vector bundle $E$ on a {\em generic} principally
polarized Abelian variety $A$ carries a Hermitian--Einstein connection
$a_E$ that defines a holomorphic structure on $E$ with curvature 
 \begin{equation}
 F_{a_E}=\La i\cdot\om,
 \label{eq7.1}
 \end{equation}
where $\La$ is any constant element of $\U(\rk(E))$, for example, $\id$.
Hence the restriction of $a_E$ to every fibre of $\pi$ is a flat Hermitian
connection on a $g$-torus, and thus
 \begin{equation}
 \begin{gathered}
 (a_E)\rest{\pi^{-1}(b)}=\chi_1\oplus\dots\oplus \chi_{\rk E},
 \quad\text{where}\\
 \chi_i\in\Hom(\pi_1(\pi^{-1}(b)),\U(1)))=\Pic(\pi^{-1}(b))=(\pi')^{-1}(b).
 \end{gathered}
 \label{eq7.2}
 \end{equation}
 \begin{dfn} For vector bundles of higher rank, an {\em $E$-{\rm BS} fibre}
is a fibre $\pi^{-1}(b)$ such that the restriction
$(E,a_E)\rest{\pi^{-1}(b)}$ admits a covariant constant section.
 \end{dfn}

Suppose that $E$ is {\em ample}, and in particular that
 \[
 H^i(A,E)=0\quad\text{for}\quad i>0.
 \]
Then, alongside the complex ``space of wave functions'' $H^0(A,E)$,
we get a new space of wave functions 
 \begin{equation}
 \sH_{\pi}^E=\bigoplus_{\text{$E$-BS}}\C\cdot s_{i},
 \label{eq7.3}
 \end{equation}
where $s_i$ is a covariant constant section of the restriction of $E$ to a
$E$-BS fibre.

We again have the problem of comparing the spaces
 \begin{equation} 
 \sH_{E}=H^0(A, E) \quad \text{and} \quad \sH_{\pi}^E.
 \label{eq7.4}
 \end{equation}

This problem can be solved by analogous (but more sophisticated) methods
from GQP. In particular

 \begin{prop} The space $H^0(A,E)$ admits a canonical theta basis.
 \end{prop}

Of course, if $X$ is any K\"ahler manifold with some real polarization
(\ref{eq1.3}) and stable holomorphic vector bundle $E$ admitting an
Hermitian connection with curvature of the form (\ref{eq7.1}), we get two
spaces
 \begin{equation} 
 \sH_{E}=H^0(X, E) \quad \text{and} \quad \sH_{\pi}^E
 \label{eq7.5}
 \end{equation}
to compare. 

In particular if $X=R_{\Si_g}$ one has
 \begin{prop} For a stable vector bundle of higher rank $E$ on $R_{\Si}$,
 \[
 H^0(R_{\Si_g}, E)=\bigoplus_{E\text{-}\BS}\C\cdot s_{i}.
 \]
 \end{prop}
This holds in particular for all the symmetric powers of the Poincar\'e
bundles.

\subsection*{K3 surfaces} If the transcendental lattice $T_S$ of a K3
surface $S$ contains an even unimodular sublattice $H$ of rank 2, then $S$
admits a real polarization (see for example \cite{G1}, \cite{G2} or
\cite{T}). Then every ample complete linear system on $S$ admits a special
theta basis.

\subsection*{Geometry behind these constructions: mirror reflection of
holomorphic geometry} If our real polarization (\ref{eq1.3}) is regular,
that is, the differential of the map $\pi\colon X\to B$ is surjective then
all fibres are $n$-tori, and the second fibration $\pi'\colon X'\to B$ can
be defined fibrewise in the usual way:
 \[
 (\pi')^{-1}(b)=\Hom(\pi_1(\pi^{-1}(b)),\U(1))
 \quad\text{for any point $b\in B$;}
 \]
that is, the fibre of $\pi'$ is the space of classes of flat connections
on the trivial line bundle on $\pi^{-1}(b)$. This is a fibration in
groups, and we want to consider its zero section $s\colon B\to X'$ as a
submanifold $s_0\subset X'$.

The restriction of a pair $(E, a_E)$ to any fibre $\pi^{-1}(b)$ defines a
finite set of points $(\pi')^{-1}(b)$, and hence a multisection
 \begin{equation}
 \GFT(E)\subset X',
 \label{eq8.1}
 \end{equation} 
which we again consider as a middle dimensional submanifold of $X'$. This
cycle is called the {\em Geometric Fourier Transformation} of $E$.

 Under the identification $s_0=B$, the set of $E$-BS fibres is defined now
as the set of intersection points
 \begin{equation}
E\text{-BS}=s_0 \cap \GFT(E),
 \label{eq8.2}
 \end{equation} 
and under some geometric conditions, we expect that the number
 \begin{equation}
 \# E\text{-BS}=[s_0] \cap [\GFT(E)]
 \label{eq8.3}
 \end{equation} 
where $[\ ]$ is the cohomology class of a submanifold.

 In the general case of a polarization with degenerate fibres, this
construction can be performed over the open subset $B_0\subset B$ of
smooth tori and a number of questions arise: 
 \begin{enumerate}
 \item to construct a smooth compactification $S'$;
 \item to construct a symplectic form $\om'$ and extend it to $S'$, in such
a way that $\pi'$ is a new real polarization;
 \item to construct a complex polarization of $S'$ such that the fibration
$\pi$ is given by construction we have described, starting from
$(S',\om',L',a')$.
 \end{enumerate}
In full generality these problems are very hard (see for example
\cite{G1}, \cite{G2}).

The ideal picture is described by the mirror diagram
 \begin{equation}
 \renewcommand{\arraystretch}{1.5}
 \begin{matrix}
 &&S &\longleftarrow & E \\
 && \kern3mm \big\downarrow \pi &&\\
 & & B&& \\
 && \kern4.5mm \big\uparrow \pi'\\
 \GFT(E)&\longrightarrow &S' &\longleftarrow & s_0 \\
 \end{matrix}
 \label{eq8.4}
 \end{equation}
with holomorphic objects (vector bundles) corresponding to the top of
(\ref{eq8.4}) and special Lagrangian cycles to the bottom.

\subsection*{The inverse problem}
Every stable holomorphic vector bundle $E$ on an $S_I$ (top of
(\ref{eq8.4})) is a point in the moduli space of stable holomorphic vector
bundles
 \begin{equation}
 E\in \sM_{[E]}^\mathrm{s}
 \label{eq8.5}
 \end{equation}
of topological type $[E]$.

But the cycle $\GFT(E)$ (bottom of (\ref{eq8.4})) is a point in the moduli
space of special Lagrangian cycles (see \cite{HL})
 \begin{equation}
 \GFT(E)\in \sM^{[\GFT(E)]}
 \label{eq8.6}
 \end{equation}
of topological type $[\GFT(E)]$. Thus we get a map
 \begin{equation}
 \GFT \colon \sM_{[E]}^\mathrm{s}\to \sM^{[\GFT(E)]}
 \label{eq8.7}
 \end{equation}
sending $E$ to $\GFT(E)$.

However, we have not used all the information contained in $E$. Namely,
$\GFT(E)$ can be defined as a {\em supercycle} (or {\em brane}). It's easy
to see that any cycle $\GFT(E)$ admits a tautological topologically trivial
line subbundle $L$ with Hermitian connection $s_{\tau}$. A pair
 \begin{equation}
(\GFT(E),a_{\tau})=\sGFT(E)
 \label{eq8.8}
 \end{equation}
of this type is called a {\em supercycle} (or brane). 

The attempt to reconstruct the vector bundle $E$ (top of (\ref{eq8.4}))
from the supercycle $\sGFT(E)$ (bottom of (\ref{eq8.4})) is called the {\em
inverse problem}. More formally, let $S\sM^{[\GFT(E)]}$ be the moduli
space of supercycles of topological type $[\GFT(E)]$. Then in many special
cases, one can prove that the map
 \[
 \sGFT\colon\sM^\mathrm{s}_{[E]}\to S\sM^{([\GFT(E)])}
 \]
is an embedding at the general point. That is, a general stable vector
bundle $E$ (top of (\ref{eq8.4})) is uniquely determined by the supercycle
$\sGFT(E)$ on $S'$ (bottom of (\ref{eq8.4})).

For example, if the fibration $\pi \colon X\to B$ is the family of all
deformations (with degenerations) of the general fibre $\pi^{-1}(b)=T^n$ as
a torus with special structure inside $S$, then
 \[
 X'=S\sM^{[T^n]}
 \] 
is the family of all deformations (with degenerations) of the pair $(T^n,
\tau_0)$, where $\tau_0$ is the trivial connection. 

At present this program is only realized in part (see, for example,
\cite{T}). We must first use the experience of the geometric quantization
procedure, and apply it in the Calabi--Yau realm of {\em simply connected
K\"ahler manifolds with canonical class zero}. But in this paper, we want
to emphasize that there exists the collection of singular Fano varieties
$R_g$ for which these constructions are very important, although this is
an extremely irregular case.

\subsubsection*{Acknowledgments}

I would like to express my gratitude to the Institut de Math\'ematiques de
Jussieu and the Ecole Normale Sup\'erieure, and personally to Joseph Le
Potier and Arnaud Beauville for support and hospitality. I wish to thank
Yves Laslo and Christoph Sorger for many helpful discussions. Thanks are
again due to Miles Reid for tidying up the English.

\bigskip
\noindent
Andrei Tyurin, Algebra Section, Steklov Math Institute,\\
Ul.\ Gubkina 8, Moscow, GSP--1, 117966, Russia \\
e-mail: Tyurin@tyurin.mian.su {\em or} Tyurin@Maths.Warwick.Ac.UK\\
{\em or} Tyurin@mpim-bonn.mpg.de
 \end{document}